\documentclass[12pt]{article}
\usepackage{amssymb,amsmath,amsopn,amsfonts}
\usepackage[dvips]{color}
\usepackage{colordvi,multicol}
\usepackage{epsf}
\newfam\msbfam
\font\tenmsb=msbm10 \textfont\msbfam=\tenmsb \font\sevenmsb=msbm7
\scriptfont\msbfam=\sevenmsb \font\fivemsb=msbm5
\scriptscriptfont\msbfam=\fivemsb

\def\th#1{\vspace{1mm}\noindent{\bf #1}\quad}

\def\proof{\vspace{1mm}\noindent{\it Proof}\quad}

\numberwithin{equation}{section}

\def\bc{\begin{center}}
\def\ec{\end{center}}
\def\no{\noindent}
\def\hang{\hangindent\parindent}
\def\textindent#1{\indent\llap{\qquad #1\ \ \enspace}\ignorespaces}
\def\ref{\par\hang\textindent}

\begin{document}

\title{ {\bf A note on stochastic semilinear equations and their associated Fokker-Planck equations
\thanks{Research supported in part by NSFC (No. 11301026) and China Postdoctoral Science Foundation funded project(2012M520153) and  the DFG through IRTG 1132 and CRC 701}\\} }
\author{ {\bf Michael R\"{o}ckner}$^{\mbox{c}}$, {\bf Rongchan Zhu}$^{\mbox{a,c},}$, {\bf Xiangchan Zhu}$^{\mbox{b,c},}$
\date{}
\thanks{E-mail address: roeckner@math.uni-bielefeld.de(M. R\"{o}ckner),
zhurongchan@126.com(R. C. Zhu), zhuxiangchan@126.com(X. C. Zhu)}\\\\
$^{\mbox{a}}$Department of Mathematics, Beijing Institute of Technology, Beijing 100081,  China \\
$^{\mbox{b}}$School of Sciences, Beijing Jiaotong University, Beijing 100044, China
 \\
$^{\mbox{c}}$Department of Mathematics, University of Bielefeld, D-33615 Bielefeld, Germany}

\maketitle
\begin{abstract}
The main purpose of this paper is to prove existence and uniqueness of (probabilistically weak and strong) solutions to stochastic differential equations (SDE) on Hilbert spaces under a new approximation condition on the drift, recently proposed in [BDR10] to solve Fokker-Planck equations (FPE), extended in this paper to a considerably larger class of drifts. As a consequence we prove existence of martingale solutions to the SDE (whose time marginals then solve the corresponding FPE). Applications include stochastic semilinear partial differential equations with white noise and a non-linear drift part which is the sum of a Burgers-type part and a reaction diffusion part. The main novelty is that the latter is no longer assumed to be of at most linear, but of at most polynomial growth. This case so far had not been covered by the existing literature. We also give a direct and more analytic proof for existence of solutions to the corresponding FPE, extending the technique from [BDR10] to our more general framework, which in turn requires to work on a suitable Gelfand triple rather than just the Hilbert state space.

\end{abstract}

\vspace{1mm}
\no{\footnotesize{\bf 2000 Mathematics Subject Classification AMS}:\hspace{2mm} 60H15,60J60,47D07}
 \vspace{2mm}

\no{\footnotesize{\bf Keywords}:\hspace{2mm}  Fokker-Planck equations, stochastic PDEs, Kolmogorov operators, martingale solutions}

\section{Introduction}

 Let $H$ be a separable real Hilbert space with inner product
$\langle\cdot, \cdot\rangle$ and corresponding norm $|\cdot|$. $L(H)$ denotes the set of all bounded linear operators on $H$, $\mathcal{B}(H)$
its Borel $\sigma$-algebra.

  Consider the following type of non-autonomous stochastic differential equations on $H$ and time interval
$[0,T]$:
\begin{equation}\left\{\begin{array}{l}dX(t)=(AX(t)+F(t,X(t)))dt+\sqrt{G}dW(t),\\
 X(s)=x\in H, t\geq s.\end{array}\right.\end{equation}
Here $W(t), t\geq0$, is a cylindrical Wiener process on $H$ defined on a stochastic basis $(\Omega, \mathcal{F},\{\mathcal{F}_t\}_{t\geq0},$ $P)$,
$G$ is a linear symmetric positive definite operator in $H$, $D(F)\in \mathcal{B}([0,T]\times H)$,
  $F:D(F)\subset [0,T]\times H\rightarrow
H$ is a Borel measurable map, and $A:D(A)\subset H\rightarrow H$ is the infinitesimal generator of a $C_0$-semigroup $e^{tA}, t\geq0$, on $H$.

Even in this case, where the noise is additive, it is a fundamental question in the theory of stochastic differential equations (SDE) in infinite dimensional state spaces with numerous applications to concrete (non-linear) stochastic partial differential equations (SPDE), whether there exists a (unique) weak or strong (in the probabilistic sense) solution to SDE (1.1). In [BDR10], for a large class of semigroup generators $A$ and in the fully elliptic case, i.e., where $G$ has an inverse $G^{-1}\in L(H)$ (in particular including the case of space-time white noise), a quite general approximation condition on $F$ was identified, which implies that (at least) the corresponding Fokker-Planck equation (FPE) has a solution, which is also unique under some $L^2$-conditions on $F$ (see [BDR11] and the recent paper [BDRS13]). The purpose of this paper is to generalize this result under the same approximation condition on $F$ (see Hypothesis 2.2 (i), (ii)) in two ways:

\textbf{(a)} We prove that (1.1) has indeed a weak (=martingale) solution (in the sense of Stroock-Varadhan). In particular, its time marginals solve the corresponding FPE.

\textbf{(b)} We prove \textbf{(a)} in a more general framework, namely allowing (as is usual in the variational approach to SDE on Hilbert spaces, see e.g. [G98], [GR00] and also [PR07]) that $F$ takes values in a larger space, more precisely in $D((-A)^{1/2})^*$, assuming (as in [BDR10]) that $A$ is negative definite and self-adjoint. In short: we shall work in a Gelfand triple.

This is done in Section 2 of this paper and the corresponding main result is summarized in Theorem 2.3 there. In order to include degenerate cases, where e.g. $\rm{Tr}G<\infty$, we assume, instead of the requirement $G^{-1}\in L(H)$ from [BDR10], that the approximating equations have martingale solutions (see Hypothesis 2.2 (iii) below), which can be easily checked in many applications.

We would like to stress, however, that, though \textbf{(b)} above may hint in this direction, our result is not at all covered by the variational approach to SDE on Hilbert spaces (see e.g. [PR07]), since, first, there is no monotonicity condition on $F$ and, second, the noise coefficient operator $G$ is not assumed to have finite trace.

We would also like to stress that in our main application (see Section 4 below) by a standard result from the seminal paper [MR99] we can also prove uniqueness of the martingale solutions. This, however, by principle cannot generally imply uniqueness of solutions to the corresponding FPE, because the latter might have solutions which are not the time marginals of a martingale solution. However, it is well-known that uniqueness for FPE implies uniqueness of martingale problems (see [SV79]). Therefore, as FPE are concerned our uniqueness results in this paper are much weaker and, in fact, far from those in [BDR11] and [BDRS13] for FPE.

Our more general framework, indicated under \textbf{(b)} above, considerably widens the range of applications in comparison with those in [BDR10].

Let us briefly describe a class of examples, which we present in detail in Section 4 of this paper and which have been studied intensively in the literature, however, under more stringent assumptions (on the function $f$ in (1.2) below).

Consider the stochastic semilinear partial differential equation (SPDE)
\begin{equation}dX(t)=(\frac{\partial^2}{\partial \xi ^2}X(t)+f(t,X(t))+\frac{\partial}{\partial \xi }g(t,X(t )))dt+\sqrt{G}dW(t ),\end{equation}
on $H:=L^2(0,1)$ with Dirichlet boundary condition
$$X(t,0)=X(t,1)=0, t\in[0,T],$$
and  initial condition
$$X(0)=x\in H,$$
where $f(\xi,t,z), g(\xi,t,z)$ are Borel measurable functions of $(\xi,t,z)\in[0,1]\times\mathbb{R}_+\times \mathbb{R}$,
$W$ is a cylindrical Wiener process on $H$ and $G$ is a linear symmetric positive definite operator in $H$.

 This kind of stochastic partial differential equations has been studied intensively. If $f=0$ and $g=\frac{1}{2}r^2$, the above equation is
 just the  stochastic Burgers equation and has been investigated in many papers (see e.g. [DDT94], [DZ96] and the references therein). When $g=0$
 then the above equation is a stochastic reaction-diffusion equation which has also attracted a lot of attention (see e.g. [DZ92], [D04],
 [BDR10] and the references therein). In the general case, this kind of equations has been studied e.g. in [G98], [GR00], where, however,
 $f$ was assumed to be of at most linear growth. We stress that the linear growth of $f$ cannot be dropped  in [G98], [GR00], since the approximation technique used there requires  this assumption.

As an application of our main result (Theorem 2.3 below) we obtain that (1.2) has a martingale solution which under a natural integrability condition is unique (see Theorem 4.2), where we assume the usual conditions on the "Burgers-part" $g$ of the drift, but in contrast to [G98], [G00] we can allow $f$ to be of polynomial growth. We, however, pay a price for considering such more general $f$, because we do not recover all results from [G98], [G00] where e.g. (1.2) with multiplicative noise is included under certain assumptions, $g$ is allowed to be of polynomial growth in [G00] and under local Lipschitz assumptions on $f$ (and $g$) also existence and uniqueness of strong solutions is shown.
If, however, we assume one sided local Lipschitz  assumption on $f$  (see (4.14) below),  we also get existence and uniqueness of strong solutions under only polynomial growth conditions on $f$ (see Theorem 4.7 below) by proving pathwise uniqueness and applying the Yamada-Watanable Theorem. We also stress that our condition for $f$ is more general than the one imposed in the corresponding applications in [BDR10] (see condition (f2) in Section 4), which allows us to take more general $f$ (see Example 4.0).

At least if $\rm{Tr} G<\infty$, we can also apply our framework to a lot of other stochastic semilinear equations, as e.g. the stochastic 2D Navier-Stokes equation (see Remark 4.9). Since in this case there are many known existence results (cf. [GRZ09] and the references therein ) based on It\^{o}'s formula and the Burkholder-Davis-Gundy inequality  to obtain the estimates required for tightness of the distributions of the approximations, we do not give details here, but concentrate on (1.2) in our applications.

Though, as mentioned above, our Theorem 2.3 implies the existence of solutions to the FPE associated to (1.1), we nevertheless also give an alternative direct proof for the latter which is more analytic in nature and a generalization of the corresponding one in [BDR10]. We think that this proof is of sufficient independent interest. Therefore,  we include it here, stressing the modifications needed in our (in comparison with that in [BDR10]) more general framework.

We mention here that recently, there has  been quite an interest in Fokker-Planck
equations with irregular coefficients in finite dimensions (see e.g. [A04], [DPL89], [F08], [BDR08a] and the references therein). In [BDR08b],
[BDR09] and [BDR10], Bogachev, Da Prato and the authors  have started the  study of Fokker-Planck equations also in infinite
dimensions, more precisely, on Hilbert spaces. Let us briefly present the formulation of the FPE corresponding to (1.1) in our framework.

The
 Kolmogorov operator $L_0$ corresponding to (1.1) reads as follows:
$$L_0u(t,x):=D_tu(t,x)+\frac{1}{2}\rm{Tr}[GD^2u(t,x)]+\langle x, A^*Du(t,x)\rangle+{ }_{V^*}\!\langle F(t,x),Du(t,x)\rangle_V,\textrm{ } (t,x)\in D(F),$$
where $D_t$ denotes the derivative in time and $D, D^2$ denote the first- and second-order Frechet derivatives in space, i.e., in $x\in H$,
respectively. Furthermore, $V:=D((-A)^{1/2})$, $V^*$ is its dual and ${ }_{V^*}\!\langle \cdot,\cdot\rangle_V$ denotes their dualization, assuming again that $A$ is negative definite and self-adjoint. The operator $L_0$ is defined on the space $D(L_0):=\mathcal{E}_A([0,T]\times H)$, defined to be the linear span of all real and
imaginary parts of all functions $u_{\phi,h}$ of the form
$$u_{\phi,h}(t,x)=\phi(t)e^{i\langle x,h(t)\rangle}, t\in [0,T], x\in H,$$
where $\phi\in C^1([0,T])$, $\phi(T)=0$, $h\in C^1([0,T];D(A^*))$ and $A^*$ denotes the adjoint of $A$.

For a fixed initial time $s\in[0, T]$ the Fokker-Planck equation is an equation
for measures $\mu(dt,dx)$ on $[s, T] \times H$ of the type
$$\mu(dt, dx)=\mu_t(dx)dt,$$
with $\mu_t\in\mathcal{P}(H)$ for all $t\in[s, T]$, and $t\mapsto\mu_t(A)$ measurable on $[s, T]$
for all $A\in\mathcal{B}(H)$, i.e. $\mu_t(dx), t\in[s, T]$, is a probability kernel from
$([s, T],\mathcal{B}([s, T]))$ to $(H,\mathcal{B}(H))$. Then the FPE corresponding to (1.1) for an initial condition
$\zeta\in\mathcal{P}(H)$ reads as follows: $\forall u\in D(L_0)$
\begin{equation}\int_Hu(t,y)\mu_t(dy)=\int_Hu(s,y)\zeta(dy)+\int_s^tds'\int_HL_0u(s',y)\mu_{s'}(dy),\textrm{ for }dt-a.e.
t\in[s,T],\end{equation}
where the dt-zero set may depend on $u$.

In Section 3 of this paper, we  prove directly the existence of  solutions to FPE (1.3) within the same framework as in Section 2, which generalizes the results in [BDR10] . In Section 4 as an
application we  prove the existence of  solutions for the FPE  associated with concrete SPDE of type (1.2), i.e. allowing
polynomially growing nonlinearities for the reaction-diffusion part $f$ and Burgers type nonlinearities $g$ at the same time (see Theorem 4.3 below), which can  not
be handled within the framework of [BDR10].

Finally, we recall that our work  covers the case  $G^{-1}\in L(H)$, i.e. the case of full (including
white) noise. If $\rm{Tr} G<\infty$, there are many other known existence results for FPE (cf. [BDR08b, BDR09] ), based on the method of constructing Lyapunov
functions with weakly compact level sets for the Kolmogorov operator $L_0$. These techniques so far could, however, not be used when $\rm{Tr} G=\infty$.

\section{Existence of martingale solutions}
Let us start with formulating our  assumptions on the coefficients of SDE (1.1).
\vskip.10in
\th{Hypothesis 2.1} (i) $A$ is self-adjoint such that there exists $\omega\in (-\infty,0)$ such that
$\langle Ax,x\rangle\leq \omega |x|^2, x\in D(A),$ and $A^{-1}$ is compact on $H$.

(ii) $G\in L(H)$ is symmetric, nonnegative.

(iii) There exists $\delta,\delta_1>0$ such that $$\int_0^T\rm{Tr}[(-A)^\delta e^{rA}G(-A)^\delta e^{rA}] dr<\infty, \quad
\int_0^1r^{-2\delta_1}\rm{Tr}[ e^{rA}G e^{rA}] dr<\infty.$$
\vskip.10in

Under Hypothesis 2.1, there exists an orthonormal basis $\{e_k\}_{k\geq0}$ for $H$ consisting of eigenfunctions of $-A$ such that the associated sequence of eigenvalues $\{\lambda_k\}$ form an increasing unbounded sequence. It is well known (see [D04, Theorem 2.9]) that under Hypothesis 2.1 (iii) the stochastic convolution
$$W_A(t)=\int_0^te^{(t-r)A}\sqrt{G}dW(r),\quad t\geq0,$$
is a well-defined  continuous  process in $H$ with values in $D((-A)^{\delta})$ and
\begin{equation}\sup_{t\in[0,T]}E|(-A)^\delta W_A(t)|^2\leq \int_0^T\rm{Tr}[(-A)^\delta e^{rA}G(-A)^\delta e^{rA}] dr<\infty.\end{equation}
\vskip.10in
\th{Remark}  If $(-A)^{2\delta-1}$ is of trace-class for some $\delta>0$ and $G\in L(H)$, Hypothesis 2.1 (iii) is obviously satisfied. We would like to point out here that there is a misprint in Hypothesis 2.1 (iii) in [BDR10], where $(-A)^{-2\delta}$ should be replaced by $(-A)^{2\delta-1}$. Likewise in the right hand side of inequality (2.1) in [BDR10].
\vskip.10in

We weaken resp. modify Hypotheses 2.2, 2.3 in [BDR10] as follows: let
$V:=D((-A)^{1/2})$ and consider the following Gelfand triple:
$$D(A)\subset V\subset H\cong H^*\subset V^*\subset D(A)^*,$$
where $V^*$ and $D(A)^*$ are the dual of $V, D(A)$ respectively and ${ }_{D(A)^*}\!\langle\cdot,\cdot\rangle_{D(A)}={ }_{V^*}\!\langle\cdot,\cdot\rangle_V=\langle\cdot,\cdot\rangle$ if restricted to $H\times D(A)$. We have the following formulas for the norm in $V, V^*$,
$$|\cdot|_V^2=\sum_k\lambda_k|\langle \cdot,e_k\rangle|^2,\quad |\cdot|_{V^*}^2=\sum_k\lambda_k^{-1}|\langle \cdot,e_k\rangle|^2.$$
Furthermore, we relax the assumptions on $F$ in (1.1) to be just $V^*$-valued. More precisely, let $F: D(F)\subset [0,T]\times H\rightarrow V^*$
be Borel measurable. Then the Kolmogorov operator is given as follows
$$L_0u(t,x):=D_tu(t,x)+\frac{1}{2}\rm{Tr}[GD^2u(t,x)]+\langle x, ADu(t,x)\rangle+{ }_{V^*}\!\langle F(t,x),Du(t,x)\rangle_V,$$
for $u\in D(L_0)$. Below we fix $s\in[0,T]$ as starting time.

\vskip.10in
\th{Hypothesis 2.2} There exist measurable maps $F_\alpha: [s,T]\times H\rightarrow D(A)^*, \alpha\in (0,1]$,
 $K>0$ and a lower semicontinuous function $J:[s,t]\times H\rightarrow [1,\infty]$ , such that the following four conditions are satisfied for all $\alpha\in(0,1]$:

 (i) for all $(t,x)\in D(F)$ and all $h\in D(A)$
$$F_\alpha(t,x)\in V^*, \quad |F_\alpha(t,x)|_{V^*}\leq J(t,x)<\infty,$$
$$|{ }_{V^*}\!\langle F(t,x)-F_\alpha(t,x),h\rangle_{V}|\leq \alpha c(h)J^2(t,x),$$
for some constant $c(h)>0$.

(ii) $(t,x)\mapsto{ }_{D(A)^*}\!\langle F_\alpha(t,x), h\rangle_{D(A)} \textrm{ is continuous on } [s,T]\times H, \forall h\in D(A), \alpha\in (0,1].$

(iii) The following approximating stochastic equations for $\alpha\in(0,1]$
\begin{equation}dX_\alpha(t)=[AX_\alpha(t)+F_\alpha(t,X_\alpha(t))]dt+\sqrt{G}dW(t), X_\alpha(s)=x\in H,\end{equation}
have a martingale solution  which we denote by $X_\alpha(\cdot,s,x)$, i.e.  there exists a stochastic basis $(\Omega,\mathcal{F},\{\mathcal{F}_t\}_{t\in [s,T]},P)$, a cylindrical
Wiener process $W$ on $H$ and a progressively measurable process $X_\alpha:[s,T]\times \Omega\rightarrow H$, such that for  $P$-a.e. $\omega\in
\Omega$ and $\phi\in D(A)$,
$$X_\alpha(\cdot,\omega)\in L^2([s,T];H)\cap C([s,T];D(A)^*),$$
$$\langle X_\alpha(t)-x,\phi\rangle=\int_s^t(\langle X_\alpha(\tau),A\phi\rangle+{}_{D(A)^*}\!\langle F_\alpha(\tau,X_\alpha(\tau)),\phi\rangle_{D(A)} ) d\tau+\int_s^t\langle\phi,\sqrt{G}dW(\tau)\rangle,\quad \forall t\in[s,T].$$

(iv)$|F|_{V^*}\leq J$ on $[s,T]\times H$, where we set $|F|_{V^*}:=+\infty$ on $[s,T]\times H\backslash D(F)$, and setting
$$P_{s,t}^\alpha\varphi(x):=E[\varphi(X_\alpha(t,s,x))], \quad 0\leq s<t\leq T, \varphi\in\mathcal{B}_b(H),$$
we have
$$\int_s^tP_{s,s'}^\alpha J^2(s',\cdot)(x)ds'\leq K\int_s^tJ^2(s',x)ds', \forall x\in H, t\in [s,T],\alpha\in (0,1].$$
\vskip.10in

\th{Remark} (i) Since $J\equiv\infty$ on $[s,T]\times H\setminus D(F)$, the latter inequality obviously holds if it holds on $D(F)$. Therefore, if we can find a function which is a Lyapunov function for $P_{s,t}^\alpha$ uniformly in $\alpha$ i.e.
$$P_{s,t}^\alpha J^2(t,\cdot)(x)\leq KJ^2(t,x),\qquad\forall (t,x)\in D(F), t\in [s,T],\alpha\in (0,1],$$ Hypothesis 2.2 (iv) is satisfied.

(ii) If $G$ has a bounded inverse and if the approximation in Hypothesis 2.2  can be chosen  such that $F_\alpha$ are bounded measurable maps, then we can use Girsanov's theorem to obtain the existence of a martingale solution. For the case that $\rm{Tr}G<\infty$, we could choose $F_\alpha=P_{[\frac{1}{\alpha}]+1} F$, where $P_n$ is the orthogonal projection onto the linear space spanned by the first $n$ eigenvectors $e_k$. Then we can apply the results in [PR07, Chapter 4] to the equation
$$dX_\alpha(t)=[AX_\alpha(t)+F_\alpha(t,X_\alpha(t))]dt+GdW(t),$$ provided $F_\alpha$ satisfies the monotonicity assumptions specified there, and  obtain the existence of a martingale solution required in Hypothesis 2.2 (iii).

(iii) In Hypothesis 2.3 (iii) the stochastic basis $(\Omega,\mathcal{F},\{\mathcal{F}_t\}_{t\in [s,T]},P)$ and the cylindrical
Wiener process $W$ may depend on $\alpha$. However, this will not change our proof since we want to prove  the laws of $X_\alpha$ are tight in a suitable space.
\vskip.10in

\th{Theorem 2.3} Assume Hypotheses 2.1, 2.2.
Then for every $x\in B:=\{x\in H: \int_s^TJ^2(t,x)dt<\infty\}$, there exists  a martingale solution to (1.1), i.e.  there exists a stochastic basis $(\Omega,\mathcal{F},\{\mathcal{F}_t\}_{t\in [s,T]},P)$, a cylindrical
Wiener process $W$ on $H$ and a progressively measurable process $X:[s,T]\times \Omega\rightarrow H$, such that for  $P$-a.e. $\omega\in
\Omega$ and $\phi\in D(A)$,
$$X(\cdot,\omega)\in L^2([s,T];H)\cap C([s,T];D(A)^*),$$
and
$$\langle X(t)-x,\phi\rangle=\int_s^t(\langle X(\tau),A\phi\rangle+{}_{D(A)^*}\!\langle F(\tau,X(\tau)),\phi\rangle_{D(A)} ) d\tau+\int_s^t\langle\phi,\sqrt{G}dW(\tau)\rangle\quad \forall t\in[0,T].$$
Moreover, for $\delta_2:=\delta\wedge\frac{1}{2}$ with $\delta$ as in Hypothesis 2.1
$$E\int_s^T(J^2(s',X(s'))+|(-A)^{\delta_2} X(s')|^2+|X(s')|^2)ds'\leq C\int_s^T(J^2(s',x)+|x|^2)ds'.$$
\proof For simplicity we assume $s=0$. For $\alpha\in(0,1]$, set $X_\alpha(t):=X_\alpha(t,0,x), x\in B,$ and
$$Y_\alpha(t):=X_\alpha(t)-W_A(t), \quad t\geq 0.$$
Then for $\phi\in D(A)$, we have $$\langle Y_\alpha(t)-x,\phi\rangle=\int_0^t(\langle Y_\alpha(s'),A\phi\rangle+{}_{D(A)^*}\!\langle F_\alpha(s',X_\alpha(s')),\phi\rangle_{D(A)} ) ds'\quad \forall t\in[0,T].$$

Choosing $\phi=e_k$ in the above equation and using Newton-Leibniz formula ,  we obtain
$$\langle Y_\alpha(t),e_k\rangle^2=\langle x,e_k\rangle^2+2\int_0^t\langle Y_\alpha(s'),e_k\rangle(\langle Y_\alpha(s'),Ae_k\rangle+{}_{D(A)^*}\!\langle F_\alpha(s',X_\alpha(s')),e_k\rangle_{D(A)} ) ds',\quad \forall t\in[0,T].$$
Then by the Cauchy-Schwarz inequality and, since $Ae_k=-\lambda_ke_k$, we have
$$\langle Y_\alpha(t),e_k\rangle^2+\int_0^t\lambda_k\langle Y_\alpha(s'),e_k\rangle^2ds'\leq\langle x,e_k\rangle^2+\int_0^t\lambda_k^{-1}|{}_{D(A)^*}\!\langle F_\alpha(s',X_\alpha(s')),e_k\rangle_{D(A)}|^2 ds'.$$
Summing over $k$ we get
$$|Y_\alpha(t)|^2+\int_0^t|(-A)^{1/2}Y_\alpha(s')|^2ds'\leq
|x|^2+\int_0^t|F_\alpha(s',X_\alpha(s'))|_{V^*}^2ds',$$
where we set $|F_\alpha|_{V^*}:=+\infty$ on $[0,T]\times H\backslash D(F)$.
Taking expectation and applying Hypothesis 2.2 yield
\begin{equation}E|Y_\alpha(t)|^2+\int_0^tE|(-A)^{1/2}Y_\alpha(s')|^2ds'\leq |x|^2+K\int_0^tJ^2(s',x)ds', \quad t\geq 0.\end{equation}
Then we deduce that for any $\varepsilon>0$ there exists $R_1>0$ such that
$$P(\int_0^T|(-A)^{1/2}Y_\alpha(s')|^2ds'>R_1)<\varepsilon, \quad \forall \alpha\in (0,1].$$
Since by Hypothesis 2.2 we have
\begin{equation}E\int_0^T|F_\alpha(s',X_\alpha(s'))|_{V^*}^2ds'\leq E\int_0^TJ^2(s',X_\alpha(s'))ds'\leq K\int_0^TJ^2(s',x)ds',\end{equation}
and
$$E\int_0^T|(-A)Y_\alpha(s')|_{V^*}^2ds'= E\int_0^T|(-A)^{1/2}Y_\alpha(s')|^2ds',$$
 we deduce that for any $\varepsilon>0$ there exists $R_2>0$ such that
$$P(\int_0^T\left|\frac{dY_\alpha}{dt}\right|_{V^*}^2ds'>R_2)<\varepsilon \quad \forall \alpha\in (0,1].$$
Then by the compactness Theorems 2.1 and 2.2 in [FG95], the laws of $X_\alpha=Y_\alpha+W_A$ are tight in $L^2([0,T];H)\cap C([0,T];D(A)^*)$.
 Thus, by Skorokhod's representation theorem there exists a subsequences $n_k$ and a sequence of random elements
$\hat{X}_k, k=1,2,3,...$ in $L^2([0,T];H)\cap C([0,T];D(A)^*)$, on some probability space $(\hat{\Omega},\hat{\mathcal{F}},\hat{P})$, such that $\hat{X}_k$ converges
almost surely in $L^2([0,T];H)\cap C([0,T];D(A)^*)$ to a random element $\hat{X}$ for $k\rightarrow\infty$ and the distributions of $\hat{X}_k$ and $X_{\frac{1}{n_k}}$
coincide. Then the second inequality in (2.4) holds for $\hat{X}_k$ and $\hat{X}$ by the lower semicontinuity of $J$.
Define for $\phi\in D(A)$,
$$\aligned\hat{M}_k(\phi)(t):=&\langle \hat{X}_k(t)-x,\phi\rangle-\int_0^t \langle\hat{X}_k(s'),A\phi\rangle ds'-\int_0^t{}_{D(A)^*}\!\langle
F_{1/n_k}(s',\hat{X}_k(s' )),\phi\rangle_{D(A)} ds'
.\endaligned$$ $\hat{M}_k(\phi)$ is a family of martingales with respect to the filtration
$$\mathcal{G}_t^k=\sigma(\hat{X}_k(r),r\leq t).$$
 For all $r\leq t\in[0,T]$ and all bounded continuous functions $\varphi$ on  $L^2([0,r];H)\cap C([0,r];D(A)^*)$ we have
$$\hat{E}((\hat{M}_k(\phi)(t)-\hat{M}_k(\phi)(r))\varphi(\hat{X}_k|_{[0,r]}))=0,$$ and $$\aligned &\hat{E}[(\hat{M}_k(\phi)(t)^2- \hat{M}_k(\phi)(r)^2-\int_r^t |\sqrt{G}\phi|^2_Hds')\varphi(\hat{X}_k|_{[0,r]})]=0.\endaligned$$
By the Burkholder-Davis-Gundy inequality we have that for $1<p<\infty$  there exists $C_p\in (0,\infty)$ such that
\begin{equation}\sup_k \hat{E}| \hat{M}_k(\phi)(t)|^{2p}\leq C_p\hat{E}(\int_0^t|\sqrt{G}\phi|_H^2dr)^p<\infty.\end{equation}
Now we prove the following estimate: for fixed $\eta>0$
\begin{equation}\hat{E}\int_0^t|{}_{D(A)^*}\!\langle F_\eta(s',\hat{X}_k(s' ))-F_\eta(s',\hat{X}(s')),\phi\rangle_{D(A)}| ds'\rightarrow0, \quad k\rightarrow\infty.\end{equation}
Indeed, we set $G_R(t,x):={}_{D(A)^*}\!\langle F_\eta(t,x),\phi\rangle_{D(A)}\chi_R({}_{D(A)^*}\!\langle F_\eta(t,x),\phi\rangle_{D(A)})$, where $\chi_R\in C_0^\infty:\mathbb{R}\rightarrow[0,1]$ is a cutoff function with
$\chi_R(r)=1$ when $|r|\leq R$ and $\chi_R(r)=0$ when $|r|>2R$. Then by the dominated convergence theorem we obtain
$$\lim_{k\rightarrow\infty}\hat{E}\int_0^t|G_R(s',\hat{X}_k(s' ))-G_R(s',\hat{X}(s'))| ds'=0.$$
Then we have $$\aligned &\lim_{R\rightarrow\infty}\sup_k\hat{E}\int_0^t|{}_{D(A)^*}\!\langle F_\eta(s',\hat{X}_k(s' )),\phi\rangle_{D(A)}-G_R(s',\hat{X}_k(s'))| ds'\\&\leq2 \lim_{R\rightarrow\infty}\sup_k\hat{E}\int_0^t|{}_{D(A)^*}\!\langle F_\eta(s',\hat{X}_k(s' )),\phi\rangle_{D(A)}|1_{\{ |{}_{D(A)^*}\!\langle F_\eta(s',\hat{X}_k(s' )),\phi\rangle_{D(A)}|>R\}}ds'\\&\leq C\lim_{R\rightarrow\infty}\sup_k\hat{E}\int_0^tJ^2(s',\hat{X}_k(s' ))ds'/R=0,\endaligned$$
where we used Hypothesis 2.2 in the second inequality and (2.4) to deduce the last convergence. The above convergence also holds for $\hat{X}$. Combining the above  estimates (2.6) follows.

By Hypothesis 2.2
we have\begin{equation}\aligned&\hat{E}\int_0^t|{}_{D(A)^*}\!\langle F_{1/n_k}(s',\hat{X}_k(s' ))-F(s',\hat{X}(s')),\phi\rangle_{D(A)}| ds'\\\leq&\hat{E}\int_0^t|{}_{D(A)^*}\!\langle F_{1/n_k}(s',\hat{X}_k(s' ))-F(s',\hat{X}_k(s')),\phi\rangle_{D(A)}| ds'\\&+\hat{E}\int_0^t|{}_{D(A)^*}\!\langle F(s',\hat{X}_k(s' ))-F_\eta(s',\hat{X}_k(s')),\phi\rangle_{D(A)}| ds'\\&+\hat{E}\int_0^t|{}_{D(A)^*}\!\langle F(s',\hat{X}(s' ))-F_\eta(s',\hat{X}(s')),\phi\rangle_{D(A)}| ds'\\&+\hat{E}\int_0^t|{}_{D(A)^*}\!\langle F_\eta(s',\hat{X}_k(s' ))-F_\eta(s',\hat{X}(s')),\phi\rangle_{D(A)}| ds'
\\\leq&C\hat{E}\int_0^t\frac{1}{n_k}J^2(s',\hat{X}_k(s')) ds'+C\hat{E}\int_0^t\eta J^2(s',\hat{X}_k(s')) ds'+C\hat{E}\int_0^t(\eta J^2(s',\hat{X}(s')) ds'\\&+\hat{E}\int_0^t|{}_{D(A)^*}\!\langle F_\eta(s',\hat{X}_k(s' ))-F_\eta(s',\hat{X}(s')),\phi\rangle_{D(A)}| ds'\\\rightarrow&0, k\rightarrow\infty,\endaligned\end{equation}
where in the second inequality we use Hypothesis 2.2 and the last convergence follows by (2.4) for $\hat{X}_k$ and $\hat{X}$ and (2.6). In fact, we could choose $\eta_0$ small enough such that the second term and the third term in the right hand side of last inequality converge to $0$. Then for such $\eta_0$ we could find $k$ large enough such that the first term and the last term converge to $0$.
Then by (2.5) and (2.7) we obtain
$$\lim_{k\rightarrow\infty}\hat{E}|\hat{M}_k(\phi)(t)-M(\phi)(t)|=0$$
and$$\lim_{k\rightarrow\infty}\hat{E}|\hat{ M}_k(\phi)(t)-M(\phi)(t)|^2=0,$$
where $$\aligned M(\phi)(t):=&\langle \hat{X}(t)-x,\phi\rangle-\int_0^t \langle\hat{X}(s'),A\phi\rangle ds'-\int_0^t\langle
F(s,\hat{X}(s' )),\phi\rangle ds'.\endaligned$$
Taking the limit we obtain that for all $r\leq t\in[0,T]$ and all bounded continuous functions $\varphi$ on   $L^2([0,r];H)\cap C([0,r];D(A)^*)$
 $$\hat{E}((M(\phi)(t)-M(\phi)(r))\varphi(\hat{X}\upharpoonright_{[0,r]}))=0.$$
and
$$\hat{E}((M(\phi)(t)^2- M(\phi)(r)^2-\int_r^t| \sqrt{G}v|_H^2 ds)\varphi(\hat{X}\upharpoonright_{[0,r]}))=0.$$
Thus the existence of a martingale solution for (1.1) follows by a martingale representation theorem (cf. [DZ92, Theorem 8.2],[O05, Theorem 2]). The last inequality follows by (2.1), (2.3), (2.4) and the lower semicontinuity of $J^2+|(-A)^{\delta_2}\cdot|^2+|\cdot|^2$.
  $\hfill\Box$

  Set
$$P_{s,t}\varphi(x):=E[\varphi(X(t,s,x))], \quad 0\leq s<t\leq T, \varphi\in\mathcal{B}_b(H),$$
   and
$$\mu_t(dx):=(P_{s,t})^*\zeta(dx),$$
where $\zeta\in\mathcal{P}(H)$ such that
$$\int_s^T\int_H(J^2(s',x)+|x|^2)\zeta(dx)ds'<\infty.$$

Now  It\^{o}'s formula implies that this is a solution to the corresponding Fokker-Planck equation, i.e.
$\forall u\in D(L_0)$
$$\int_Hu(t,y)\mu_t(dy)=\int_Hu(s,y)\zeta(dy)+\int_s^tds'\int_HL_0 u(s',y)\mu_{s'}(dy),\textrm{ for all }
t\in[s,T].$$

\section{Existence of solutions to the Fokker-Planck equation}

In this section we prove directly the existence of solutions for the Fokker-Planck equation (1.3) under the same conditions as in Section 2.

Set $$W_A(t,s)=\int_s^te^{(t-s')A}\sqrt{G}dW(s'), \quad t\geq s.$$

The Kolmogorov operator $L_\alpha$ corresponding to (2.2) is given by
$$\aligned L_\alpha u(t,x):=&D_tu(t,x)+\frac{1}{2}\rm{Tr}[GD^2u(t,x)]+\langle x, ADu(t,x)\rangle\\&+{}_{D(A)^*}\!\langle F_\alpha (t,x),Du(t,x)\rangle_{D(A)},\textrm{ } (t,x)\in [0,T]\times H,\quad u\in D(L_0).\endaligned$$

Fix $s\in[0,T)$ and set
$$\mu_t^\alpha(dx):=(P_{s,t}^\alpha)^*\zeta(dx),$$
where $\zeta\in\mathcal{P}(H)$ is the initial condition, at $t=s$.

Now  It\^{o}'s formula implies that this is a solution to the corresponding Fokker-Planck equation, i.e.
$\forall u\in D(L_0)$
\begin{equation}\int_Hu(t,y)\mu_t^\alpha(dy)=\int_Hu(s,y)\zeta(dy)+\int_s^tds'\int_HL_\alpha u(s',y)\mu_{s'}^\alpha(dy),\textrm{ for all }
t\in[s,T],\end{equation}

\vskip.10in

\th{Theorem 3.1} Assume Hypotheses 2.1, 2.2 and let $\zeta\in \mathcal{P}(H)$ be such that
$$\int_s^T\int_H(J^2(s',x)+|x|^2)\zeta(dx)ds'<\infty.$$ Then there exists a solution $\mu_t(dx)dt$ to the Fokker-Planck equation (1.3) such that
$$\sup_{t\in[s,T]}\int_H|x|^2\mu_t(dx)<\infty,$$
and
$$t\mapsto\int_Hu(t,x)\mu_t(dx)$$
is continuous on $[s,T]$ for all $u\in D(L_0)$. Finally, for some $C>0$ and for $\delta_2:=\delta\wedge\frac{1}{2}$ with $\delta$ as in Hypothesis 2.1 one has
\begin{equation}\int_s^T\int_H(J^2(s',x)+|(-A)^{\delta_2} x|^2+|x|^2)\mu_{s'}(dx)ds'\leq C\int_s^T\int_H(J^2(s',x)+|x|^2)\zeta(dx)ds'.\end{equation}
\proof For $\alpha\in(0,1]$, set $X_\alpha(t):=X_\alpha(s,t,x), x\in H,$ and
$$Y_\alpha(t):=X_\alpha(t)-W_A(t,s), \quad t\geq s.$$
By the same arguments to obtain  (2.3) we also have here that
\begin{equation}E|Y_\alpha(t)|^2+\int_s^tE|(-A)^{1/2}Y_\alpha(s')|^2ds'\leq |x|^2+K\int_s^tJ^2(s',x)ds', \quad t\geq s.\end{equation}
Then  for $s\leq t\leq T$ we obtain
$$E|X_\alpha(t)|^2\leq 2|x|^2+2K\int_s^TJ^2(s',x)ds'+2\kappa,$$
where $\kappa:=\sup_{t\in[s,T]}E|W_A(t)|^2<\infty$. Now we integrate with respect to $\zeta$ over $x\in H$ and obtain for all $s\leq t\leq T$ and some $C\in (0,\infty)$ that
\begin{equation}\int_H|x|^2\mu_t^\alpha(dx)\leq C[1+\int_s^T\int_H(J^2(s',x)+|x|^2)\zeta(dx)ds'].\end{equation}
Hence we can use Prohorov' theorem  (see [B07, Theorem 8.6.7]) to obtain that  for each $t\in[s,T],$ there exists a
sub-sequence $\{\alpha_n\}$ (possibly depending on $t$) such that the measures $\mu_t^{\alpha_n}$ converge $\tau_w$-weakly to a measure $\tilde{\mu}_t\in\mathcal{P}(H)$ as
$n\rightarrow\infty$, where $\tau_w$ denotes the weak topology on $H$.

Now we have that for $\varphi\in\mathcal{E}_A(H)$, defined to be the set of all linear combinations of all real parts of functions of the form $x\mapsto e^{i\langle x,h\rangle}, h\in D(A)$,
\begin{equation}t\mapsto \mu_t^\alpha(\varphi):=\int_H\varphi(x)\mu_t^\alpha(dx), \alpha\in(0,1]\textrm{ are equicontinuous on }[s,T].\end{equation} In fact, for $s\leq t_1\leq t_2\leq T$
$$\aligned|\mu_{t_2}^\alpha(\varphi)-\mu_{t_1}^\alpha(\varphi)|\leq& \frac{1}{2}\|\rm{Tr}[GD^2\varphi]\|_\infty|t_2-t_1|\\
&+|t_2-t_1|^{1/2}\|AD\varphi\|_\infty(\int_{t_1}^{t_2}\int_H|x|^2\mu_{s'}^\alpha(dx)ds')^{1/2}\\
&+|t_2-t_1|^{1/2}\|(-A)^{1/2}D\varphi\|_\infty(\int_{t_1}^{t_2}\int_HJ^2(s',x)\zeta(dx)ds')^{1/2},\endaligned$$
where $\|\cdot\|_\infty$ denotes the sup-norm on $H$. By (3.4) and Hypothesis 2.2, (3.5) follows.

Then by the same arguments as in the proof of [BDR10, Theorem 2.6], we can construct a measure $\mu_t$ and a subsequence $\{\alpha_n\}$ such that
$\mu_t^{\alpha_n}$ converge $\tau_w$-weakly to $\mu_t$ for all $t\in[0,T]$. Indeed, by a diagonal argument we can choose $\{\alpha_n\}$ such that $\mu_t^{\alpha_n}\rightarrow\tilde{\mu}_t$ $\tau_\omega$-weakly as $n\rightarrow\infty$ for every rational $t\in[s,T]$. Moreover (3.4) holds for $\tilde{\mu}_t$ in place of $\mu_t^\alpha$ for $t\in[s,T]\cap\mathbb{Q}$. Hence by [B07, Theorem 8.6.7], for each $t\in[s,T]\setminus\mathbb{Q}$ there exists $r_n(t)\in[s,T]\cap \mathbb{Q}, n\in\mathbb{N}$ converging to $t$ and $\mu_t\in\mathcal{P}(H)$ such that $\tilde{\mu}_{r_n(t)}\rightarrow\mu_t$ $\tau_w$-weakly as $n\rightarrow\infty$. Now for fix $t\in [s,T]\backslash\mathbb{Q}$ suppose $\{\mu_t^{\alpha_n}\}$ does not weakly converge to $\mu_t$. Then by (3.4) and [B07, Theorem 8.6.7] there exists a subsequence $\{\alpha_{n_k}\}$, $\varphi\in\mathcal{E}_A(H)$ and $\nu\in\mathcal{P}(H)$ such that $\mu_t^{\alpha_{n_k}}\rightarrow\nu$ $\tau_w$-weakly as $k\rightarrow\infty$ and $\mu_t(\varphi)\neq\nu(\varphi)$. On the other hand,
for all $n,k\in\mathbb{N}$
$$\aligned|\nu(\varphi)-\mu_t(\varphi)|\leq&|\nu(\varphi)-\mu_t^{\alpha_{n_k}}(\varphi)|+\sup_{l\in\mathbb{N}}
|\mu_t^{\alpha_{n_l}}(\varphi)-\mu_{r_n(t)}^{\alpha_{n_l}}(\varphi)|\\&
+|\mu_{r_n(t)}^{\alpha_{n_k}}(\varphi)-\tilde{\mu}_{r_n(t)}(\varphi)|+|\tilde{\mu}_{r_n(t)}(\varphi)-\mu_t(\varphi)|.\endaligned$$
Letting $k\rightarrow\infty$ and then $n\rightarrow\infty$ it follows from (3.5) that $\mu_t(\varphi)=\nu(\varphi)$. Letting $\mu_t:=\tilde{\mu}_t$ for $t\in[s,T]\cap \mathbb{Q}$, we have that $\mu_t^{\alpha_n}$ converge $\tau_w$-weakly to $\mu_t$ for all $t\in[0,T]$.

 (3.4) and Lebesgue's dominated convergence theorem imply that $t\mapsto\int_Hu(t,x)\mu_t(dx)$
is continuous on $[s,T]$ for all $u\in D(L_0)$.

Now for $\delta_2:=\delta\wedge\frac{1}{2}$ with $\delta$ as in Hypothesis 2.1 (iii) by (3.3) and (2.1) we obtain
\begin{equation}\int_s^T\int_H|(-A)^{\delta_2} x|^2\mu_t^\alpha(dx)dt\leq C[1+\int_s^T\int_H(J(s',x)^2+|x|^2)\zeta(dx)ds'],\end{equation}
which implies that $\mu_t^{\alpha_n}(dx)dt$ converge weakly to $\mu_t(dx)dt$ on
$[s,T]\times H$ by the compactness of $(-A)^{-\delta_2}$. Now (3.2) follows from (3.4), (3.6) and the lower semicontinuity of $J^2+|(-A)^{\delta_2}\cdot|^2+|\cdot|^2$.

It remains to prove that $\mu_t(dx)dt$ solves the Fokker-Planck equation (1.3). Since every $h\in C^1([0,T];D(A))$ can be written as a uniform limit of piecewise affine $h_n\in C([0,T];D(A)), n\in\mathbb{N}$, it  follows by (3.2) and linearity that $\mu_t(dx)dt$ satisfies the Fokker-Planck equation (1.3) if and only if it does so for all $u\in D(L_0)$ such that $u(t,x)=\phi(t)e^{i\langle h(t),x\rangle}, x\in H,t\in[0,T]$, with $\phi\in C^1([0,T]), \phi(T)=0$ and piecewise affine $h\in C([0,T];D(A))$. Fix such a function $u\in D(L_0)$, by (3.1) we have
$$\int_s^T\int_HL_{\alpha_n}u(t,x)\mu_t^{\alpha_n}(dx)dt=-\int_s^Tu(s,x)\zeta(dx),$$
with $\alpha_n$ as above.

Since we already know that $\mu_t^{\alpha_n}(dx)dt\rightarrow\mu_t(dx)dt$ weakly and since the coefficient of the second order part of $L_{\alpha_n}$ is just $G$ (hence independent of $n$), it now suffices to prove that for all $g\in
C_b([s,T]\times H)$ and all piecewise affine $h\in C([0,T];D(A))$,
\begin{equation}
\lim_{n\rightarrow\infty}\int_s^T\int_HF_{\alpha_n}^h(t,x)g(t,x)\mu_t^{\alpha_n}(dx)dt=
\int_s^T\int_HF^h(t,x)g(t,x)\mu_t(dx)dt,\end{equation}
where
$$F_\alpha^h(t,x):={ }_{D(A)^*}\!\langle F_\alpha(t,x), h(t)\rangle_{D(A)}
+\frac{\langle Ah(t),x\rangle}{1+\alpha|\langle Ah(t),x\rangle|},$$
$$F^h(t,x):={ }_{D(A)^*}\!\langle F(t,x), h(t)\rangle_{D(A)}
+\langle Ah(t),x\rangle.$$
For $\eta\in (0,1]$ we have
\begin{equation}\aligned&|\int_s^T\int_HF_{\alpha_n}^h(t,x)g(t,x)\mu_t^{\alpha_n}(dx)dt-
\int_s^T\int_HF^h(t,x)g(t,x)\mu_t(dx)dt|\\\leq&\|g\|_\infty\int_s^T\int_H|F_{\alpha_n}^h(t,x)-F^h(t,x)|\mu_t^{\alpha_n}(dx)dt\\
&+\|g\|_\infty\int_s^T\int_H|F^h(t,x)-F_\eta^h(t,x)|\mu_t^{\alpha_n}(dx)dt\\&
+\|g\|_\infty\int_s^T\int_H|F^h(t,x)-F_\eta^h(t,x)|\mu_t(dx)dt\\&+
|\int_s^T\int_HF_\eta^h(t,x)g(t,x)\mu_t^{\alpha_n}(dx)dt-\int_s^T\int_HF_\eta^h(t,x)g(t,x)\mu_t(dx)dt|.\endaligned\end{equation}

By Hypothesis 2.2 we have for all $\alpha,\beta\in(0,1]$
\begin{equation}\aligned\int_s^T\int_H|F_\beta^h(t,x)-F^h(t,x)|\mu_t^\alpha(dx)dt\leq& \beta C(h)\int_s^T\int_H(J^2(t,x)+|x|^2)\mu_t^\alpha(dx)dt\\\leq&
 \beta C(h)C(1+\int_s^T\int_H(J^2(t,x)+|x|^2)\zeta(dx)dt),\endaligned\end{equation}
 where $C$ is a constant independent of $\alpha, \beta$ and we used Hypothesis 2.2 and (3.4) in the last step. This implies that if $n\rightarrow\infty$ and $\eta\rightarrow0$ the first two terms in (3.8) converge to zero. Since (3.9) holds for $\mu_t$ in place of $\mu_t^\alpha$, we deduce that the third term converges to zero if $\eta\rightarrow0$.
Now we consider the last summand. Since $F_\eta^h$ is continuous on $[s,T]\times H$ by our assumption, there exists a continuous function $\tilde{G}_R$ on $[s,T]\times {H}$ satisfying $\|\tilde{G}_R\|_\infty\leq R$, and $\tilde{G}_R(t,x)=F^h_{\eta}(t,x)$ on $ B_R$, for $B_R:=\{|F^h_\eta|\leq R\}$. By the weak convergence we obtain
$$\lim_{n\rightarrow\infty}\int_s^T\int_{H}\tilde{G}_R(t,x)g(t,x)\mu_t^{\alpha_n}(dx)dt=\int_s^T\int_{H}\tilde{G}_R(t,x)g(t,x)\mu_t(dx)dt.$$
By the above estimate we get
$$\aligned &\int_s^T\int_{{H}}|\tilde{G}_R(t,x)-F_{\eta}^h(t,x)|\mu_t^{\alpha_n}(dx)dt\\ \leq &CR\int_{B_R^c}\mu_t^{\alpha_n}(dx)dt+CC(h)\int_{B_R^c}(|F_\eta(t,x)|_{V^*}+|x|)\mu_t^{\alpha_n}(dx)dt\\ \leq &CR^{-1}\int_s^T\int_H(J^2(t,x)+|x|^2)\mu_t^{\alpha_n}(dx)dt+C\gamma(h)\int_{B_R^c}(J(t,x)+|x|)\mu_t^{\alpha_n}(dx)dt,\endaligned$$
where in the last inequality we used Hypothesis 2.2. Then the last summand converges to zero if $R\rightarrow\infty$ and $n\rightarrow\infty$. Hence
 (3.7) is verified and the assertion follows.
$\hfill\Box$

\section{Application}
Consider the stochastic semilinear partial differential equation
\begin{equation}dX(t)=(\frac{\partial^2}{\partial \xi ^2}X(t)+f(t,X(t))+\frac{\partial}{\partial \xi }g(t ,X(t)))dt+\sqrt{G}dW(t
),\end{equation}
with Dirichlet boundary condition
\begin{equation}X(t,0)=X(t,1)=0, t\in[0,T],\end{equation}
and the initial condition
\begin{equation}X(0)=x,\end{equation}
on $H=L^2(0,1):=L^2((0,1),d\xi)$, with $d\xi=$ Lebesgue measure. Here $f,g:(0,1)\times [0,T]\times \mathbb{R}\rightarrow\mathbb{R}$ are functions such that for every $\xi\in(0,1)$ the maps
$f(\xi,\cdot,\cdot), g(\xi,\cdot,\cdot)$ are continuous on $[0,T]\times \mathbb{R}$ and satisfy the following conditions:

\no (f1) There exist $m\in\mathbb{N}$ (without loss of generality $m\geq2$) and a nonnegative function $c_1\in L^2(0,T)$ such that for all $t\in [0,T], z\in \mathbb{R},\xi\in (0,1)$
$$|f(\xi,t,z)|\leq c_1(t)(1+|z|^m).$$
(f2) There exists a nonnegative function $c_2\in L^1(0,T)$ and $m_1\in (0,\infty)$ such that for all $t\in[0,T], z_1,z_2\in\mathbb{R}, \xi\in (0,1)$
$$(f(\xi,t,z_1+z_2)-f(\xi,t,z_1))z_2\leq c_2(t)(|z_2|^2+|z_1|^{m_1}+1).$$
(g1) The function $g$ is of the form $g(\xi,t,z)=g_1(\xi,t,z)+g_2(t,z)$, where $g_1$ and $g_2$ are Borel functions of $(\xi,t,z)\in (0,1)\times
[0,T]\times \mathbb{R}$ and of $(t,z)\in  [0,T]\times \mathbb{R}$, respectively. The function $g_1$ satisfies a linear growth and the function
$g_2$ a quadratic growth condition, i.e. there is a constant $K$ such that
$$|g_1(\xi,t,z)|\leq K(1+|z|), \quad |g_2(t,z)|\leq K(1+|z|^2),$$
for all $t\in [0,T],\xi \in (0,1), z\in \mathbb{R}$.

\no (g2) $g$ is a locally Lipschitz function with linearly growing Lipschitz constant, i.e. there exists a constant $L$ such that
$$|g(\xi,t,z_1)-g(\xi,t,z_2)|\leq L(1+|z_1|+|z_2|)|z_1-z_2|,$$
for all $t\in [0,T],\xi\in (0,1), z_1,z_2\in \mathbb{R}$.
\vskip.10in

\th{Example 4.0} Now we give examples for $f$ satisfying (f1) (f2). Let
$f:(0,1)\times [0,T]\times \mathbb{R}\rightarrow\mathbb{R}$ be a function such that for every $\xi\in(0,1)$ the maps
$f(\xi,\cdot,\cdot)$ are continuous on $[0,T]\times \mathbb{R}$. Moreover
$f=f_1+f_2$ satisfies the polynomial growth condition (f1) for some $m\geq 2$ and there exists a constant $C$ such that
$$f_1(\xi,t,\cdot)\in C^1(\mathbb{R}), \quad \partial_zf_1(\xi,t,z)\leq C, \quad (\xi,t,z)\in(0,1)\times
[0,T]\times \mathbb{R},$$
$$f_2(\xi,t,z)z\leq C[1+|z|^2], \quad |f_2(\xi,t,z)|\leq C(1+|z|^{2-\frac{1}{m}})\quad (\xi,t,z)\in(0,1)\times
[0,T]\times \mathbb{R}.$$

It immediately follows from the mean value theorem that $f_1$ satisfies (f2). Now we check (f2) for $f_2$: for $t\in[0,T], z_1,z_2\in\mathbb{R}, \xi\in (0,1)$
$$\aligned (f_2(\xi,t,z_1+z_2)-f_2(\xi,t,z_1))z_2\leq& f_2(\xi,t,z_1+z_2)(z_1+z_2)-f_2(\xi,t,z_1+z_2)z_1+(1+|z_1|^2)|z_2|\\
\leq& C+C(z_1+z_2)^2+C(1+|z_1+z_2|^{2-\frac{1}{m}})|z_1|+(1+|z_1|^m)|z_2|\\\leq&|z_2|^2+C(1+|z_1|^{2m}).\endaligned$$

\vskip.10in
Let $A:D(A)\subset H\rightarrow H$ be defined by
$$Ax(\xi)=\frac{\partial^2}{\partial\xi^2}x(\xi),\xi\in(0,1), \quad D(A)=H^2(0,1)\cap H_0^1(0,1).$$
Then $V=H_0^1(0,1)$. Let $D(F):=[0,T]\times L^{2m}(0,1)$ and for $(t,x)\in D(F)$
$$F:=F_1+F_2,\quad F_1(t,x)(\xi):=f(\xi,t,x(\xi)),\quad F_2(t,x)(\xi):={\partial_\xi g(\xi,t,x(\xi))}, \xi\in (0,1),$$
where $F_2$ takes values in $V^*$.

Finally, let $G\in L(H)$ be symmetric, nonnegative and such that $G^{-1}\in L(H)$ and there exist $\theta, q\geq0$ with $\frac{1}{2q}+2\theta<1$ such that
$$\|(\sum_k(A^{-\theta}\sqrt{G}(e_k))^2)^{1/2}\|_{L^q}<\infty,\eqno(G.1)$$
where $\{e_k\}$ is an orthonormal basis of $H$.

If $G=Id$, (G.1) is obviously satisfied. By (G.1), [B97, Corollary 3.5] and [D04, Exercise 2.16] we know that  $W_A$ is a Gaussian random variable in $C([0,T]\times [0,1])$.

 It is easily checked that $A, G$ satisfy Hypothesis 2.1 with $\delta, \delta_1\in (0,\frac{1}{4})$.

For $\alpha\in (0,1]$ and $(t,x)\in [0,T]\times H$ we define $F_\alpha:[0,T]\times H\rightarrow D(A)^*$,
 $$F_\alpha:=F_1^\alpha+F_2, \quad F_1^\alpha(t,x)(\xi):=\frac{F_1(t,x)(\xi)}{1+\alpha|F_1(t,x)(\xi)|}, \xi\in[0,1]. $$
If $F_1\equiv0$, there exists a unique (probabilistically) strong solution to (4.1) by [G98, Theorem 2.1].
  Since $F_1^\alpha$ is bounded, by  Girsanov's Theorem (cf. [MR99, Theorem 3.1], [DFPR12, Theorem 13]) , we obtain that there exists  a stochastic basis $(\Omega,\mathcal{F},\{\mathcal{F}_t\}_{t\in [0,T]},P)$, a cylindrical
Wiener process $W$ on $H$ and a progressively measurable process $X_\alpha:[s,T]\times \Omega\rightarrow H$ as in Hypothesis 2.2 (iii) satisfying  the following stochastic
differential equation
\begin{equation}dX_\alpha(t)=[AX_\alpha(t)+F_\alpha(t,X_\alpha(t))]dt+\sqrt{G}dW(t), X_\alpha(s)=x,s\leq t,\end{equation}
for all $x\in H$.

Define for $m\geq 2$ as in (f1)
$$J(t,x):=\left\{\begin{array}{ll}2(c_1(t)+K)(1+|x|_{L^{2m}(0,1)}^m),&\ \ \ \ \textrm{ if } (t,x)\in D(F)\\+\infty,&\ \ \ \ \textrm{
otherwise. }\end{array}\right.$$
By (g1) we have
$$|F_2(t,x)|_{V^*}\leq 2K(1+|x|_{L^4}^2)\leq J(t,x)<\infty \quad\forall (t,x)\in D(F)=[0,T]\times L^{2m}(0,1).$$
By (f1)  we obtain
$$|F(t,x)|_{V^*}\leq J(t,x)<\infty \quad\forall (t,x)\in D(F)=[0,T]\times L^{2m}(0,1).$$
One also easily checks that $F_\alpha$ satisfies Hypothesis 2.2 (i)-(iii). It remains to check the last part of Hypothesis 2.2 (iv), which, however, immediately follows from the following proposition.

\vskip.10in
\th{Proposition 4.1} For any $s\in[0,T)$, there exists $C\in(0,\infty)$, such that for $\alpha\in (0,1], x\in L^{2m}(0,1)$
$$E(|X_\alpha(t,s,x)|_{L^{2m}(0,1)}^{2m})\leq C(1+|x|_{L^{2m}(0,1)}^{2m}),\quad \forall t\in[s,T].$$
\proof Set $Y_\alpha(t):=X_\alpha(t,s,x)-W_A(s,t), t\in [s,T].$  Then we obtain for $\phi\in D(A), t\in[s,T]$
\begin{equation}\langle Y_\alpha(t)-x,\phi\rangle=\int_s^t(\langle Y_\alpha(s'),A\phi\rangle+{}_{D(A)^*}\!\langle F_\alpha(s',X_\alpha(s')),\phi\rangle_{D(A)} ) ds'.\end{equation}
 Since the trajectories of $W_A$ can be approximated by functions $W_A^n:=(1-\frac{1}{n}A)^{-1}W_A$ from $C([s,T],H^2)$ in $L^2([s,T],H)\cap C([s,T],H^{-2})$, we can replace $W_A$  by  smooth functions
$W_A^n$. Moreover, we can approximate $g$ by  smooth functions $g_n:=\varphi_n*\chi_n(g)$ for  smooth functions $\varphi_n$ on $[0,1]\times \mathbb{R}$ with supp$\varphi_n\subset [-\frac{1}{n},\frac{1}{n}]^2$ and $\chi_n:\mathbb{R}\rightarrow[0,n]$ is a smooth function on $\mathbb{R}$ satisfying $\chi_n(r)=r$ if $|r|\leq n$, $\chi_n(r)=0$ if $|r|>2n$ and $|\chi_n'|\leq C$ for a constant $C$ independent of $n$. We also approximate $x$ by smooth functions $x_n$ such that $|x_n|_{L^{2m}}\leq |x|_{L^{2m}}$.  Then each $g_n$ has bounded derivative with respect to $\xi$ and $z$ and satisfies (g1), (g2) with $K, L$ replaced by $2K,3CL$ respectively.  By a standard method ( see e.g. [GRZ09, Theorem 4.6]) we obtain that there exists a stochastic basis $(\Omega,\mathcal{F},\{\mathcal{F}_t\},P)$ and a pair process $(Y_\alpha^n, \bar{W}_A^n)$ such that $$Y_\alpha^n\in L^\infty([s,T],H)\cap L^2([s,T],H_0^1)\cap C([s,T],H^{-2})\quad P-a.s.$$ and $\bar{W}_A^n$ has the same distribution as $W_A^n$ and
for $\phi\in H_0^1, t\in[s,T]$
\begin{equation}\langle Y_\alpha^n(t)-x_n,\phi\rangle=\int_s^t({}_{V^*}\!\langle A\phi, Y_\alpha^n(s')\rangle_{V}+{}_{V^*}\!\langle F_1^\alpha(s',Y_\alpha^n(s')+\bar{W}_A^n(s'))+F_2^n(s',Y_\alpha^n(s')+\bar{W}_A^n(s')),\phi\rangle_{V} ) ds',\end{equation}
where $F_2^n(t,x)(\xi):=\partial_\xi g_n(\xi,t,x(\xi))$. Below we denote $W_A^n$ as $\bar{W}_A^n$ if there's no confusion. Now taking $\phi=\lambda_k e_k$ and $e_k$ as in (4.6)  and by the product rule for $\lambda_k\langle Y_\alpha^n(t),e_k\rangle$ and $\langle Y_\alpha^n(t),e_k\rangle$ we have
$$\aligned&\lambda_k\langle Y_\alpha^n(t),e_k\rangle^2+\int_s^t\lambda_k^2\langle Y_\alpha^n(s'),e_k\rangle^2ds'\\\leq&\lambda_k\langle x_n,e_k\rangle^2+\int_s^t|{}_{V^*}\!\langle F_1^\alpha(s',Y_\alpha^n(s')+W_A^n(s'))+F_2^n(s',Y_\alpha^n(s')+W_A^n(s')),e_k\rangle_{V}|^2 ds'.\endaligned$$
Then taking sum  we have the following  estimate since $g_n$ has bounded derivative,
$$\aligned&|Y_\alpha^{n}(t)|_V^2+\int_s^t| AY_\alpha^{n}(s')|^2ds'\\\leq &|x_n|_V^2+ \int_s^tC_\alpha+C_n(1+|Y_\alpha^{n}(s')|_V^2+|W_A^n(s')|_V^2) ds',\endaligned$$
which combining with Gronwall's lemma implies that $Y_\alpha^n\in L^\infty([s,T],H^1_0)\cap L^2([s,T],H^2)$.

Moreover, (4.6) can be easily extended to $\phi\in \{u\in L^2([0,T],H_0^1): \frac{du}{dt}\in L^2([0,T],V^*) \}$:
$$\aligned\langle Y_\alpha^n(t),\phi(t)\rangle-\langle x_n,\phi(s)\rangle=&\int_s^t({}_{V^*}\!\langle \frac{d\phi}{dt}(s'),Y_\alpha^n(s')\rangle_{V}+{}_{V^*}\!\langle A\phi(s'), Y_\alpha^n(s')\rangle_{V}\\&+{}_{V^*}\!\langle F_1^\alpha(s',Y_\alpha^n(s')+W_A^n(s'))+F_2^n(s',Y_\alpha^n(s')+W_A^n(s')),\phi(s')\rangle_{V} ) ds'\endaligned$$
Since $Y_\alpha^n\in L^\infty([s,T],H^1_0)\cap L^2([s,T],H^2)$, we can choose $\phi=(Y^n_\alpha(t))^{2m-1}$ and obtain for $t\in [s,T]$
$$\aligned &\frac{1}{2m}\frac{d}{dt}\int |Y_\alpha^n(t)|^{2m}d\xi+(2m-1)\int
|Y_\alpha^n(t)|^{2m-2}|\partial_\xi Y_\alpha^n(t)|^2d\xi\\=&\int F_1^\alpha(t,Y_\alpha^n(t)+W_A^n(s,t))
Y^n_\alpha(t)^{2m-1}d\xi+{}_{V^*}\!\langle F_2^n(t,Y^n_\alpha(t)+W^n_A(s,t)),Y^n_\alpha(t)^{2m-1}\rangle_V\\:=&I_1+I_2.\endaligned$$
Let us estimate $I_2$. We have
\begin{equation}\aligned &{}_{V^*}\!\langle F_2^n(t,Y^n_\alpha(t)+W^n_A(s,t)),Y^n_\alpha(t)^{2m-1}\rangle_V\\=&{}_{V^*}\!\langle
[F_2^n(t,Y^n_\alpha(t)+
W^n_A(s,t))-F_2^n(t,Y^n_\alpha(t))],Y^n_\alpha(t)^{2m-1}\rangle_V+{}_{V^*}\!\langle  F_2^n(t,Y^n_\alpha(t)),Y^n_\alpha(t)^{2m-1}\rangle_{V}
.\endaligned\end{equation}
For the first term on the right hand side of (4.7), we have by (g2), and Young's inequality
$$\aligned &{}_{V^*}\!\langle [F_2^n(t,Y^n_\alpha(t)+
W^n_A(s,t))-F_2^n(t,Y^n_\alpha(t))],Y^n_\alpha(t)^{2m-1}\rangle_V\\\leq& C\int (1+|Y^n_\alpha(t)|+|W_A^n(s,t)|)|W_A^n(s,t)||Y_\alpha^n(t)|^{2m-2}
|\partial_\xi Y_\alpha^n(t)|d\xi\\\leq&\frac{1}{2}\int|Y_\alpha^n(t)|^{2m-2}|\partial_\xi
Y_\alpha^n(t)|^2d\xi+C\int|W_A^n(s,t)|^2|Y_\alpha^n(t)|^{2m}d\xi\\&+C\int (1+|W_A^n(s,t)|)|W_A^n(s,t)|
|Y^n_\alpha(t)|^{2m-2}|\partial_\xi Y^n_\alpha(t)|d\xi\\\leq&\int |Y^n_\alpha(t)|^{2m-2}|\partial_\xi
Y^n_\alpha(t)|^2d\xi+C|W_A^n(s,t)|_{L^{4m}}^{4m}+C|W_A^n(s,t)|_{L^{2m}}^{2m}+
(C|W_A^n(s,t)|_{L^\infty}^2+C)|Y_\alpha^n(t)|_{L^{2m}}^{2m}.\endaligned$$
For the second term on the right hand side of (4.7), we have
$$\aligned\int_0^1 g_2^n(t,Y^n_\alpha) Y^n_\alpha(t)^{2m-2}\partial_\xi Y^n_\alpha(t)d\xi=\int_0^1 \partial_\xi g_3(t,Y^n_\alpha)
d\xi=0,\endaligned$$
where $g_3(t,r)=\int_0^rg_2^n(t,z)z^{2m-2}dz$. Then we obtain by (g1)
$$\aligned{}_{V^*}\!\langle  F_2^n(t,Y^n_\alpha(t)),Y^n_\alpha(t)^{2m-1}\rangle_{V}=&-(2m-1)\int g_1^n(\xi,t,Y^n_\alpha)
Y^n_\alpha(t)^{2m-2}\partial_\xi Y^n_\alpha(t)d\xi\\\leq& C\int (1+
|Y^n_\alpha(t)|)|Y^n_\alpha(t)|^{2m-2}|\partial_\xi Y^n_\alpha(t)|d\xi\\\leq&\int (|Y^n_\alpha(t)|^{2m-2}|\partial_\xi
Y^n_\alpha(t)|^2+C|Y^n_\alpha(t)|^{2m}+C)d\xi.\endaligned$$
 Now we consider $I_1$. We note that by (f1), (f2), for all $y,z\in\mathbb{R}, t\in[0,T],\xi\in(0,1),$
 $$\aligned&|f(\xi,t,y+z)y|\\
 =&|(f(\xi,t,y+z)-f(\xi,t,z))y+f(\xi,t,z)y|\\\leq&c_2(t)(|y|^2+|z|^{m_1}+1)+c_1(t)(1+|z|^{m})|y|
 \\\leq&c(t)(1+|y|^2+|z|^{m_1}+|z|^m|y|),\endaligned$$
 where $c(t)=c_1(t)+c_2(t)$.
 Then $$I_1\leq c(t)\int [1+|Y_\alpha^n(t)|^2+|W_A^n(s,t)|^m|Y_\alpha^n(t)|+|W_A^n(s,t)|^{m_1}]|Y_\alpha^n(t)|^{2m-2}d\xi.$$
 Now we obtain
\begin{equation}\aligned &\frac{1}{2m}\frac{d}{dt}\int |Y^n_\alpha(t)|^{2m}d\xi+\int
|Y_\alpha^n(t)|^{2m-2}|\partial_\xi Y_\alpha^n(t)|^2d\xi\\\leq &
c(t)\int[1+(2+\frac{2m-1}{2m})|Y^n_\alpha(t)|^{2m}+\frac{1}{2m}|W_A^n(s,t)|^{2m^2}+\frac{1}{m}|W_A^n(s,t)|^{mm_1}]d\xi\\&+ C|W_A^n(s,t)|_{L^{4m}}^{4m}
+C|W_A^n(s,t)|_{L^{2m}}^{2m}+C+(C|W_A^n(t)|_{L^{\infty}}^2+c)|Y^n_\alpha(t)|_{L^{2m}}^{2m}
.\endaligned\end{equation}
 Then $c(\cdot)\in L^1([0,T])$ by (f1), (f2) and Gronwall's lemma yields that
\begin{equation}\aligned|Y_\alpha^n(t)|_{L^{2m}}^{2m}\leq &e^{\int_s^t
C|c(t')|+C(|W_A^n(s,t')|_{L^\infty}^2+1)dt'}(|x_n|_{L^{2m}}^{2m}+C\int_s^t(|c(t')|(|W_A^n(s,t')|_{L^{2m^2}}^{2m^2}\\&+|W_A^n(s,t')|_{L^{mm_1}}^{mm_1})+|W_A^n(s,t')|_{L^{4m}}^{4m}+
|W_A^n(s,t')|_{L^{2m}}^{2m}+1) dt')\\\leq &e^{\int_s^t
C|c(t')|+C(|W_A(s,t')|_{L^\infty}^2+1)dt'}(|x|_{L^{2m}}^{2m}+C\int_s^t(|c(t')|(|W_A(s,t')|_{L^{2m^2}}^{2m^2}\\&+|W_A(s,t')|_{L^{mm_1}}^{mm_1})+|W_A(s,t')|_{L^{4m}}^{4m}+
|W_A(s,t')|_{L^{2m}}^{2m}+1) dt').\endaligned\end{equation}
By (G.1), [B97, Corollary 3.5] and [D04, Exercise 2.16] we know that  $W_A$ is a Gaussian random variable in $C([s,T]\times [0,1])$. Then we have for any $p>1$
$$E\sup_{(t,\xi)\in[s,T]\times (0,1)}|W_A(s,t)(\xi)|^p<\infty,$$
and by Fernique's Theorem (cf. [DZ92, Theorem 2.6])  there exists a constant $\varepsilon>0$ independent of $s$ such that
\begin{equation}Ee^{\varepsilon \sup_{(t,\xi)\in[s,T]\times (0,1)}|W_A(s,t)(\xi)|^2}<\infty,\end{equation}
where $\varepsilon>0$ can be chosen independent of $s$. Indeed, since by the Markov property of $W_A$ we have for any $r>0$ that
$$P(\sup_{t\in[s,T]}|W_A(s,t)|_{L^\infty}\leq r)=P(\sup_{t\in[s,T]}|W_A(t-s)|_{L^\infty}\leq r)\geq P(\sup_{t\in[0,T]}|W_A(t)|_{L^\infty}\leq r),$$
we can choose common $\varepsilon$ and $r$ such that
$$\log(\frac{1-P(\sup_{t\in[0,T]}|W_A(t)|_{L^\infty}\leq r)}{P(\sup_{t\in[0,T]}|W_A(t)|_{L^\infty}\leq r)})+32\varepsilon r^2\leq -1.$$
Then (4.10) follows from  Fernique's Theorem.

Taking expectation in (4.9) we obtain for $s\leq t\leq t_0$ such that $t_0-s$ is small enough,
\begin{equation}E|Y_\alpha^n(t)|_{L^{2m}}^{2m}\leq C|x|_{L^{2m}}^{2m}+C,\end{equation}
where $C$ is a constant independent of $\alpha,n$. By (4.8) and (4.9) we have
$$E|Y_\alpha^n(t)|^2_{L^{2}([s,t_0],H^{1})}\leq C|x|_{L^{2m}}^{2m}+C.$$
Moreover, since by (f1) (g1)  we have
$$\aligned&\int_s^{t_0}(|AY_\alpha^n(s')|_{H^{-1}}^2+|F_1^\alpha(s',Y_\alpha^n(s')+W_A^n(s'))|_{H^{-1}}^2+|F_2^n(s',Y_\alpha^n(s')+W_A^n(s'))|_{H^{-1}}^2) ds'\\\leq& C\int_s^{t_0}|Y_\alpha^n(s')|_{H^1}^2ds'+C\int_s^{t_0}(c_1(s')^2+K^2)(1+|Y_\alpha^n(s')+W_A^n(s')|_{L^{2m}}^{2m})ds',\endaligned$$
we obtain
$$E|Y_\alpha^n(t)|^2_{W^{1,2}([s,t_0],H^{-1})}\leq C|x|_{L^{2m}}^{2m}+C.$$

Thus by [FG95, Theorems 2.1, 2.2] we get $Y_\alpha^n$ in $L^2([s,t_0],H)\cap C([s,t_0],H^{-2})$  are tight. Also $W_A^n$ in $L^2([s,t_0],H)\cap C([s,t_0],H^{-2})$ are tight. Therefore, we have $(Y_\alpha^n,W_A^n)$
are tight in $(L^2([s,t_0],H)\cap C([s,t_0], H^{-2}))\times (L^2([s,t_0],H)\cap C([s,t_0], H^{-2}))$. Hence there exists a subsequence (still denoted by $(Y^n_\alpha, W_A^n)$) converging  in distribution. By Skorohod's embedding theorem, there exist a stochastic basis $(\tilde{\Omega},\tilde{\mathcal{F}},\{\tilde{\mathcal{F}}_t\}_{t\in [s,t_0]}, \tilde{P})$ and, on this basis, $L^2([s,t_0];H)\cap C([s,t_0], H^{-2})$-valued random variables $\tilde{Y}^n_\alpha,\tilde{Y}_\alpha,\tilde{W}_A^n,\tilde{W}_A, n\geq1$, such that for every $n\in\mathbb{N}$, $(\tilde{Y}^n_\alpha,\tilde{W}_A^n)$ has the same law as $(Y^n_\alpha,W_A^n)$ on $(L^2([s,t_0];H)\cap C([s,t_0], H^{-2}))\times (L^2([s,t_0];H)\cap C([s,t_0], H^{-2}))$, and $\tilde{Y}^n_\alpha\rightarrow \tilde{Y}_\alpha, \tilde{W}_A^n\rightarrow \tilde{W}_A$ in $L^2([s,t_0];H)\cap C([s,t_0], H^{-2}),\tilde{P}$ -a.s.. Then (4.11) holds for $\tilde{Y}^n_\alpha,\tilde{Y}_\alpha$.
 For each $n\geq1$, define the process $$\aligned\tilde{M}_n(t):=&\tilde{Y}_\alpha^n(t)+\tilde{W}_A^n(t)-x_n-\int_s^tA\tilde{Y}_\alpha^n(s')ds'-
 \int_s^tA\tilde{W}_A^n(s')ds'-\int_s^t F_1^\alpha(s',\tilde{Y}_\alpha^n(s')+\tilde{W}_A^n(s'))ds'\\&-\int_s^t F_2^n(s',\tilde{Y}_\alpha^n(s')+\tilde{W}_A^n(s'))ds'.\endaligned$$
 In fact $\tilde{M}_n$ is a square integrable martingale with respect to the filtration $$\{\mathcal{G}_n\}_t=\sigma\{\tilde{Y}^n_\alpha(r), \tilde{W}_A^n(r),r\leq t\}.$$ For all $r\leq t\in[s,t_0]$, all bounded continuous functions $\phi$ on  $(C([s,r];H^{-2})\cap L^2([s,r];H))\times (C([s,r];H^{-2})\cap L^2([s,r];H))$ and all $v\in C^\infty([0,1])$, we have
$$\tilde{E}(\langle \tilde{M}_n(t)-\tilde{M}_n(r),v\rangle\phi(\tilde{Y}_\alpha^n\upharpoonright_{[s,r]},\tilde{W}_A^n\upharpoonright_{[s,r]}))=0$$ and $$\aligned &\tilde{E}((\langle \tilde{M}_n(t),v\rangle^2-\langle \tilde{M}_n(r),v\rangle^2-\int_r^t |(1-\frac{1}{n}A)^{-1}\sqrt{G}v|^2_Hds)\phi(\tilde{Y}_\alpha^n\upharpoonright_{[s,r]},\tilde{W}_A^n\upharpoonright_{[s,r]}))=0.\endaligned$$
By the Burkholder-Davis-Gundy inequality we have for $1<p<\infty$
$$\sup_n \tilde{E}|\langle \tilde{M}_n(t),v\rangle|^{2p}\leq C\sup_nE(\int_0^t|(1-\frac{1}{n}A)^{-1}\sqrt{G}v|_H^2dr)^p<\infty.$$Since $\tilde{Y}_\alpha^n\rightarrow \tilde{Y}_\alpha, \tilde{W}_A^n\rightarrow \tilde{W}_A$ in $L^2([s,t_0];H)\cap C([s,t_0], H^{-\beta})$, we have
$$\aligned&\tilde{E}\int_s^t |\langle F_2^n(s',\tilde{Y}_\alpha^n(s')+\tilde{W}_A^n(s'))-F_2(s',\tilde{Y}_\alpha(s')+\tilde{W}_A(s')),v\rangle| ds'\\\leq&\tilde{E}\int_s^t |\langle F_2^n(s',\tilde{Y}_\alpha^n(s')+\tilde{W}_A^n(s'))-F_2^n(s',\tilde{Y}_\alpha(s')+\tilde{W}_A(s')),v\rangle| ds'\\&+\tilde{E}\int_s^t |\langle F_2^n(s',\tilde{Y}_\alpha(s')+\tilde{W}_A(s'))-F_2(s',\tilde{Y}_\alpha(s')+\tilde{W}_A(s')),v\rangle| ds'\\\leq&C\tilde{E}\int_s^t(|\tilde{Y}_\alpha^n(s')|+|\tilde{W}_A^n(s')|+|\tilde{W}_A(s')|+1+|\tilde{Y}_\alpha(s')|)[|\tilde{Y}_\alpha^n(s')-\tilde{Y}_\alpha(s')|+
|\tilde{W}_A^n(s')-\tilde{W}_A(s')|]ds'\\&+\tilde{E}\int_s^t |\langle F_2^n(s',\tilde{Y}_\alpha(s')+\tilde{W}_A(s'))-F_2(s',\tilde{Y}_\alpha(s')+\tilde{W}_A(s')),v\rangle| ds'\\\rightarrow&0, \textrm{ as }n\rightarrow\infty,\endaligned$$
where in the second inequality we used (g2) and we used (4.11) to obtain the  convergence.
The other terms can be estimated similarly, which altogether implies that
$$\lim_{n\rightarrow\infty}\tilde{E}|\langle\tilde{ M}_n(t)-M(t),v\rangle|=0$$
and$$\lim_{n\rightarrow\infty}\tilde{E}|\langle\tilde{ M}_n(t)-M(t),v\rangle|^2=0,$$
where $$M(t):=\tilde{Y}_\alpha(t)+\tilde{W}_A(t)-x-\int_s^tA\tilde{Y}_\alpha(s')ds'-
 \int_s^tA\tilde{W}_A(s')ds'-\int_s^t F_\alpha(s',\tilde{Y}_\alpha(s')+\tilde{W}_A(s'))ds.$$
Taking the limit we obtain that for all $r\leq t\in[s,t_0]$, all bounded continuous functions on  $(C([s,r];H^{-\beta})\cap L^2([s,r];H))\times (C([s,r];H^{-\beta})\cap L^2([s,r];H))$ and  $v\in C^\infty([0,1])$,
 $$\tilde{E}(\langle M(t)-M(r),v\rangle\phi(\tilde{Y}_\alpha\upharpoonright_{[s,r]},\tilde{W}_A\upharpoonright_{[s,r]}))=0.$$
and
$$\tilde{E}((\langle M(t),v\rangle^2-\langle M(r),v\rangle^2-\int_r^t| \sqrt{G}v|_H^2 ds')\phi(\tilde{Y}_\alpha\upharpoonright_{[s,r]},\tilde{W}_A\upharpoonright_{[s,r]}))=0.$$
Thus, the existence of a martingale solution for (4.4) follows by a martingale representation theorem (cf. [DZ92, Theorem 8.2],[O05, Theorem 2]).
Now we obtain $\tilde{X}_\alpha=\tilde{Y}_\alpha+\tilde{W}_A$ is a martingale solution of (4.4) in $[s,t_0]$. Thus, by Girsanov's Theorem and the pathwise uniqueness of the solution to (4.4) when $f\equiv0$, we obtain the  uniqueness of (the distributions for) the martingale solution of
(4.4), which implies  that $\tilde{X}_\alpha$ has the same distribution as $X_\alpha$.
By this and (4.11) we have for $s\leq t\leq t_0$,
$$E|Y_\alpha(t)|_{L^{2m}(0,1)}^{2m}\leq C|x|_{L^{2m}(0,1)}^{2m}+C.$$
Moreover,
\begin{equation}E|X_\alpha(t, s,x)|_{L^{2m}(0,1)}^{2m}\leq C|x|_{L^{2m}(0,1)}^{2m}+C,\end{equation}
where $C$ is a constant independent of $\alpha,s$.
Furthermore, by [EK86, Theorem 4.2]  and the  uniqueness of the distributions for the martingale solution of (4.4) we obtain that  the laws of the martingale solutions $X_\alpha(t,s,x)$ of (4.4) form a Markov process. We use
$\mu^\alpha_{s,t}(x,dy)$ to denote the distribution of $X_\alpha(t,s,x)$, $x\in H$. Then by the Markov property we have for $0\leq s\leq t_1\leq t_2\leq T, x\in H$
$$\mu^\alpha_{s,t_2}(x,dz)=\int_H \mu^\alpha_{s,t_1}(x,dy)\mu^\alpha_{t_1,t_2}(y,dz).$$
By this and (4.12) we obtain by iteration that for any $t\in[s,T]$
$$\int|z|_{L^{2m}}^{2m}\mu^\alpha_{s,t}(x,dz)=\int\int |z|_{L^{2m}}^{2m} \mu^\alpha_{t_1,t}(y,dz)\mu^\alpha_{s,t_1}(x,dy)\leq
C|x|_{L^{2m}(0,1)}^{2m}+C,$$
which is exactly our assertion.
$\hfill\Box$
\vskip.10in

 Since $c_1\in L^2([0,T])$ by (f1), the set $B$ in Theorem 2.3 is $L^{2m}(0,1)$. By Theorem 2.3 we now obtain the following:
\vskip.10in
\th{Theorem 4.2}Suppose that (f1), (f2), (g1), (g2), (G.1) hold. For each initial value $x\in L^{2m}(0,1)$ there exists a martingale
solution to problem (4.1)-(4.3), i.e.  there exists a stochastic basis $(\Omega,\mathcal{F},\{\mathcal{F}_t\}_{t\in [0,T]},P)$, a cylindrical
Wiener process $W$ on $H$ and a progressively measurable process $X:[0,T]\times \Omega\rightarrow H$, such that for  $P$-a.e. $\omega\in
\Omega$,
$$X(\cdot,\omega)\in L^\infty([0,T];L^{2}(0,1))\cap C([0,T];H^{-2})$$ and  for all $\phi\in C^2([0,1])$
$$\aligned \langle X(t),\phi\rangle =&\langle x, \phi\rangle +\int_0^t\langle X(r),\partial^2_\xi\phi\rangle dr+\int_0^t\langle f(r ,X(r)),\phi\rangle dr
\\&-\int_0^t\langle g(r ,X(r )),\partial_\xi\phi\rangle dr+\int_0^t \langle\phi,  \sqrt{G}dW(r )\rangle \quad\forall t\in[0,T]\quad
P-a.s..\endaligned$$
We also have \begin{equation}X-W_A\in L^2([0,T],H^1)\quad P-a.s.,\quad E\int_0^T|X(t)|_{L^{2m}}^{2m}dt<\infty.\end{equation}
Moreover, if $P, P'$ are the laws of two martingale solutions on $C([0,T];H^{-2})$ to problem (4.1)-(4.3) with the same initial value $x\in L^{2m}$
and  $$\int_0^T |\omega(t)|_{L^{2m}}^{2m}dt<\infty \quad P+P'-a.s.,$$
then $P=P'$, where $\omega(\cdot)$ is the canonical process on $C([0,T];H^{-2})$.

\proof (4.13) follows from (4.8)-(4.11).  The weak uniqueness follows by (f1),  [MR99, Theorem 3.3] and the pathwise uniqueness of the solution of (4.4) when $f\equiv0$. Here we can extend a solution to $C([0,\infty),H^{-2})$ by taking $X(t)=X(T), t\geq T$ and apply the results in [MR99]. $\hfill\Box$
\vskip.10in

Likewise, Theorem 3.1 applies to  all $\zeta\in \mathcal{P}(H)$  such that
$$\int_H|x|^{2m}_{L^{2m}(0,1)}\zeta(dx)<\infty.$$
More precisely, we have:

\th{Theorem 4.3} Let $\zeta\in \mathcal{P}(H)$ be such that
$$\int_H|x|_{L^{2m}}^{2m}\zeta(dx)<\infty.$$ Then there exists a solution $\mu_t(dx)dt$ to the Fokker-Planck equation (1.3) such that
$$\sup_{t\in[s,T]}\int_H|x|^2\mu_t(dx)<\infty$$
and
$$t\mapsto\int_Hu(t,x)\mu_t(dx)$$
is continuous on $[s,T]$ for all $u\in D(L_0)$. Finally, for some $C>0$ and $\delta\in(0,\frac{1}{4})$ as in Hypothesis 2.1 (iii)
$$\int_s^T\int_H(|x|_{L^{2m}}^{2m}+|(-A)^\delta x|^2+|x|^2)\mu_{r}(dx)dr\leq C\int_H|x|_{L^{2m}}^{2m}\zeta(dx).$$

\th{Remark 4.4} (i) Here we choose the $L^{2m}$-norm as a Lyapunov function $J$ in Hypothesis 2.2. In [RS06], the first named author of this paper and
Sobol studied  the above semilinear stochastic partial differential equations with time independent coefficients. They also choose the $L^{2m}$-norm as a Lyapunov
function with weakly compact level sets for the Kolmogorov operator $L_0$ and by analyzing the resolvent of the operator $L$ they constructed
a unique martingale solution to this problem if the noise is trace-class. In this paper, we concentrate on  space-time white noise for which  the method
of constructing Lyapunov functions with weakly compact level sets for the Kolmogorov operator $L_0$ is more delicate than in the case, where $\rm{Tr}
G<\infty$.

(ii) If $g\equiv0$, we can  obtain  the uniqueness of the solution to the Fokker-Planck equation by [BDR11, Theorem 4.1].

\vskip.10in

To obtain pathwise uniqueness, we additionally assume that  $f$ satisfies the following inequality: for $t\in[0,T], \xi\in[0,1], z_1,z_2\in\mathbb{R}$, \begin{equation}\langle f(\xi, t ,z_1)-f(\xi,t ,z_2),z_1-z_2\rangle\leq L(1+|z_1|^{m-1}+|z_2|^{m-1})|z_1-z_2|^2.\end{equation}

Now we give the definition of a (probabilistically) strong solution to (4.1)-(4.3).
\vskip.10in

\th{Definition 4.5}  We say that there exists a (probabilistically) strong solution to (4.1)-(4.3) over the time interval $[0,T]$ if for every probability space $(\Omega,\mathcal{F},\{\mathcal{F}_t\}_{t\in [0,T]},P)$ with an $\mathcal{F}_t$-Wiener process $W$, there exists  an $\mathcal{F}_t$-adapted process $X:[0,T]\times \Omega\rightarrow H$ such that
for $P-a.s.$ $\omega\in \Omega$
$$X(\cdot,\omega)\in L^\infty([0,T];L^{2}(0,1))\cap C([0,T];H^{-2}),$$  and for all $\phi\in C^2([0,1])$ we have  $P$-a.s.
$$\aligned \langle X(t),\phi\rangle =&\langle X_0, \phi\rangle +\int_0^t\langle X(r),\partial_\xi^2\phi\rangle dr+\int_0^t\langle f(r ,X(r)),\phi\rangle dr
\\&-\int_0^t\langle g(r ,X(r )),\partial_\xi\phi\rangle dr+\int_0^t \langle\phi,  \sqrt{G}dW(r )\rangle \quad \forall t\in[0,T].\endaligned$$

\vskip.10in

\th{Theorem 4.6} Suppose that $f$ satisfies (4.14). Then there exists at most one probabilistically strong solution to (4.1)-(4.3) such that $$\int_0^T|X(t)|_{L^{2m}}^{2m}dt<\infty,\quad P-a.s.,$$ and $$X-W_A\in L^2([0,T],H_0^1)\quad P-a.s..$$

\proof Consider two solutions $X_1,X_2$ of (4.1)-(4.3) in the interval $[0,T]$. Since $X-W_A\in L^2([0,T],H_0^1)$ $P$-a.s. and $X_1-X_2\in L^2([0,T],H_0^1)$ $P$-a.s., we have
$$\aligned \langle X_1(t)-X_2(t),\phi\rangle =&\int_0^t\langle X_1(r)-X_2(r),\partial^2_\xi\phi\rangle dr+\int_0^t\langle f(r ,X_1(r))-f(r,X_2(r)),\phi\rangle dr
\\&-\int_0^t\langle g(r ,X_1(r ))-g(r,X_2(r)),\partial_\xi\phi\rangle dr, \quad\forall t\in[0,T],\quad
P-a.s..\endaligned$$
Taking $\phi=e_k$ we obtain
$$\aligned \langle X_1(t)-X_2(t),e_k\rangle^2 =&2\int_0^t\langle X_1(r)-X_2(r),e_k\rangle[\langle X_1(r)-X_2(r),\partial^2_\xi e_k\rangle \\&+\langle f(r ,X_1(r))-f(r,X_2(r)),e_k\rangle
\\&-\langle g(r ,X_1(r ))-g(s,X_2(r)),\partial_\xi e_k\rangle] dr, \quad\forall t\in[0,T],\quad
P-a.s..\endaligned$$
Summing  over $k$, we obtain
$$\aligned &| X_1(t)-X_2(t)|^2+2\int_0^t|\nabla(X_1(r)-X_2(r))|^2dr\\\leq&2\int_0^t\langle f(r ,X_1(r))-f(r,X_2(r)),X_1(r)-X_2(r)\rangle dr
\\&-2\int_0^t\langle g(r ,X_1(r ))-g(s,X_2(r)),\partial_\xi(X_1(r)-X_2(r))\rangle dr.
\endaligned$$
For the first term on the right hand side by (4.14) we have
$$\aligned&\int_0^t\langle f(r ,X_1(r))-f(r,X_2(r)),X_1(r)-X_2(r)\rangle dr\\\leq&C\int_0^t|X_1(r)-X_2(r)|_{L^4}^2(1+|X_1(r)|_{L^{2m}}^{m-1}+|X_2(r)|_{L^{2m}}^{m-1})dr\\\leq&
C\int_0^t|X_1(r)-X_2(r)|_{H^{1/4}}^2(1+|X_1(r)|_{L^{2m}}^{m-1}+|X_2(r)|_{L^{2m}}^{m-1})dr\\\leq&
C\int_0^t|X_1(r)-X_2(r)|_{H^{1}}^{1/2}|X_1(r)-X_2(r)|^{3/2}(1+|X_1(r)|_{L^{2m}}^{m-1}+|X_2(r)|_{L^{2m}}^{m-1})dr\\\leq&
\int_0^t\varepsilon|X_1(r)-X_2(r)|_{H^{1}}^2+C|X_1(r)-X_2(r)|^{2}(1+|X_1(r)|_{L^{2m}}^{2m}+|X_2(r)|_{L^{2m}}^{2m})dr,\endaligned $$
where we used Holder's inequality in the first inequality, $H^{1/4}\subset L^4$ in the second inequality, the interpolation inequality in the third inequality and Young's inequality in the last inequality.
For the second term on the right hand side we have
$$\aligned&\int_0^t\langle g(r ,X_1(r ))-g(r,X_2(r)),\partial_\xi(X_1(r)-X_2(r))\rangle dr\\\leq&C\int_0^t|\partial_\xi(X_1(r)-X_2(r))||X_1(r)-X_2(r)|_{L^4}(1+|X_1(r)|_{L^{4}}+|X_2(r)|_{L^{4}})dr
\\\leq&C\int_0^t|\partial_\xi(X_1(r)-X_2(r))|^{5/4}|X_1(r)-X_2(r)|^{3/4}(1+|X_1(r)|_{L^{4}}+|X_2(r)|_{L^{4}})dr\\\leq&
\int_0^t\varepsilon|X_1(r)-X_2(r)|_{H^{1}}^2+C|X_1(r)-X_2(r)|^{2}(1+|X_1(r)|_{L^{2m}}^{2m}+|X_2(r)|_{L^{2m}}^{2m})dr.\endaligned $$
where we used $H^{1/4}\subset L^4$ and  the interpolation inequality in the second inequality and Young's inequality in the last inequality. Combining the above three inequalities and using Gronwall-Bellman's inequality, $X_1=X_2$ follows.
$\hfill\Box$
\vskip.10in

Combining Theorems 4.3 and  4.6 we obtain the following  existence and uniqueness result by using the Yamada-Watanabe Theorem (cf. [Ku07, Theorem 3.14]).
\vskip.10in

\th{Theorem 4.7}  Suppose that (f1), (f2), (4.14), (g1), (g2), (G.1) hold.
Then for each initial condition $X_0\in L^{2m}(0,1)$, there exists a pathwise unique probabilistically strong solution $X$ of equation (4.1) over
$[0,T]$ with  initial condition $X(0)=X_0$ such that $$\int_0^T|X(t)|_{L^{2m}}^{2m}dt<\infty\quad P-a.s.$$ and \begin{equation}X-W_A\in L^2([0,T],H^1)\quad P-a.s..\end{equation}

\vskip.10in

\th{Remark 4.8} If $c_1$ in (f1) is bounded and (4.14) is modified to the following stronger local Lipschitz condition
$$| f(\xi, t ,z_1)-f(\xi,t ,z_2)|\leq L(1+|z_1|^{m-1}+|z_2|^{m-1})|z_1-z_2|,$$
condition (4.15) can be dropped. Then we can also prove that there exists a unique  probabilistically strong solution $X\in C([0,T],L^{2m})$ by considering  mild  solutions and using similar arguments as in [G98].

\th{Remark 4.9} If $\rm{Tr}
G<\infty$, we can apply Theorems 2.3 and 3.1 to other  stochastic semilinear equations and to higher dimension. For example, we can consider the 2D stochastic Navier-Stokes equation. Let $O$ be a bounded domain in $\mathbb{R}^2$ with smooth boundary. Define
$$V:=\{v\in H_0^1(O;\mathbb{R}^2),\textrm{div} v=0 \textrm{ a.e. in } O\},$$ and $H$ to be the
closure of $V$ with respect to $L^2$-norm.
The linear operator $P_H$ (Helmhotz-Hodge projection) and A (Stokes operator with viscosity
$\nu$) are defined by
$$P_H: L^2(O,\mathbb{R}^2)\rightarrow H \textrm{ orthogonal projection };$$
$$A : H^2(O,\mathbb{R}^2)\cap V\rightarrow H: Ax = \nu P_H\Delta x.$$
The nonlinear operator $F:V\rightarrow V^*$ is defined by $F(x):=-P_H[x\cdot\nabla x]$. Then if $G$ is a trace-class symmetric non-negative operator, Hypothesis 2.1 is satisfied. For Hypothesis 2.2 we  choose $F_\alpha=P_{\frac{1}{[\alpha]+1}}F$ as in Remark (ii) before Theorem 2.3 and $J(t,x):=|x||x|_V+1$. Then by It\^{o}'s formula we know that Hypothesis 2.2 (iv) is satisfied. Consequently, we   obtain the existence of a martingale solution for the stochastic 2D Navier-Stokes equation. Of course, as said in the introduction, this result is well-known and not the best possible for the 2D Navier-Stokes equation. Therefore, we omit the details here.

\th{Acknowledgement.} We thank the referee of  an earlier quite different version of this paper, whose questions and comments led to a very much improved version of this work with considerably stronger results.

\end{document}